\documentclass{amsart}

\usepackage{amsmath, amsfonts, amsthm, amssymb, amscd}

\usepackage{manfnt}

\usepackage[all]{xy}
\usepackage{color,xcolor}

\def\trd{\textcolor{red}}

\usepackage[T1]{fontenc}

\newtheorem{theorem}{Theorem}[section]

\theoremstyle{definition}
\newtheorem{definition}[theorem]{Definition}

\newtheorem{thm}{Theorem}[section]

\theoremstyle{definition}
\newtheorem{dfn}[theorem]{Definition}

\theoremstyle{remark}
\newtheorem{rmk}[theorem]{Remark}

\numberwithin{equation}{section}

\def\lh{\hbox to 15pt{\vbox{\vskip 6pt\hrule width 6.5pt height 1pt}
\kern -4.0pt\vrule height 8pt width 1pt\hfil}} 

\def\qed{\hbox{${\vcenter{\vbox{\hrule height 0.4pt\hbox{\vrule width
0.4pt height 6pt \kern5pt\vrule width 0.4pt}\hrule height 0.4pt}}}$}}

\def\Cal{\mathcal C}

\def\C2{F(S^1,2)}
\def\C3{\overbar F(S^1,3)}
\def\C4{\overbar F(S^1,4)}

\def\overbar{\overline}

\newcommand{\da}{_{\alpha}}   
   \newcommand{\db}{_{\beta}}    
   \newcommand{\dgam}{_{\gamma}}   
   \newcommand{\dab}{_{\alpha\beta}}   
   \newcommand{\dagam}{_{\alpha\gamma}}   
   \newcommand{\dbgam}{_{\beta\gamma}}

\newcommand{\Hom}{\operatorname{Hom}}

\newcommand{\p}[1]{\cal{P}_{#1}}

\def\4mu4{\nu^{a}_{\alpha\beta}}

\def\Aoo{A_\infty}

\def\p{\vskip2ex\hskip5ex}

\def\be{\begin{eqnarray}}
\def\ee{\end{eqnarray}}

\def\Z{\mathbb Z}
\def\g{\gamma}
\def\BS{Block and Smith }
\def\RH{Riemann--Hilbert }
\def\ddo{\end{document}}

\def\shrep{representation up to homotopy/homotopy coherent representation}
\def\shreps{representations up to homotopy/homotopy coherent representations}

\begin{document}
\title{\small{Homotopy coherent    representations}}
 \author{Tim Porter}
 \address{Ynys M\^{o}n /  Anglesey, Cymru / Wales, ex-University of Bangor}
  \email{t.porter.maths@gmail.com}
\author{Jim Stasheff}
\address{University of Pennsylvania}
\email{jds@math.upenn.edu}
\markboth{Homotopy coherent representations and functors \today}{\today}

\maketitle

\begin{abstract}
 \emph{Homotopy coherence} has a considerable history, albeit also by other names. 
 For this volume highlighting symmetries, the appropriate use is \emph{homotopy coherence of representations}, at one time known as \shreps.
We provide a brief semi-historical survey providing some links that may not be common knowledge.
\end{abstract}

\tableofcontents

\section{Background/Motivation: The topological setting ...}
 \emph{Homotopy coherence} has a considerable history, albeit also by other names. 
 For this volume highlighting symmetries, the appropriate use is \emph{homotopy coherence of representations}, at one time known as \shreps.
 Here is a brief semi-historical  and idiosyncratic survey providing some links that may not be common knowledge.
  \p
In topology, the fundamental theorem of covering spaces asserted that for a `nice'\footnote{globally and locally path-connected and semi-locally simply-connected}  topological space $X$,  the functor\footnote{Way back when, functors had not yet been named !} given by sending a covering space over $X$ to the corresponding representation of its fundamental group $\pi_1(X)$ as permutations of a given discrete fiber over a point is an equivalence.
That is, the representation determines 
 a covering space for which the monodromy/holonomy
 is naturally isomorphic to the original representation.
 \p
 For fiber bundles and especially for (Hurewicz) fiber spaces, there is a much more subtle correspondence. Consider such a fiber space $F\to E\to B$.
 One version of the classification of such bundles/fibrations is by the action of the based loop space $\Omega B$ on the fiber $F$. This began with Hilton showing that there was an action $\Omega B \times F \to F$ which at the level of homotopy classes gave a representation of  $\Omega B$ on $F.$
 In fact, that action gives rise to two maps  $\Omega B \times \Omega B \times F \to F$ which are homotopic.
 \p
  However, by passing to the homotopy classes of based loops in $B,$ that information is lost and there was more to be discovered. In those days,  
  $\Omega B$ was the space of based loops which were parameterized in terms of the unit interval and hence associative only up to homotopy.
   It was Masahiro Sugawara \cite{sugawara:g,sugawara:hc} who first studied `higher homotopies',  leading to today's world of $\infty$-structures of various sorts \cite{sss,loday-vallette,lurie:HighTop2006}. Although Sugawara did not treat actions up to higher homotopies as such, he did introduce \emph{strongly homotopy multiplicative maps} of associative H-spaces. With John Moore's introduction of an associative space of based loops, the adjoint $\Omega B  \to F^F$ could be seen as such a map. Considering the action as a \emph{representation }of $\Omega B$ on $F$, we have perhaps the first example of a \emph{homotopy coherent representation}.
 \begin{rmk} An apparently very different sort of treatment of \emph{higher homotopies} occurred earlier in a strong form of Borsuk's shape theory.  This was initiated by Christie in  \cite{christie:net:1944}. He used homotopy coherence in a truncated form; his work went unnoticed until the 1970s. Borsuk's form of shape theory did not fully reflect the geometric nature of his intuition.  He approximated spaces by polyhedra, up to homotopy, but did not consider the homotopies as part of the structure, just their existence. Various people tried various approaches to give a stronger form, some aspects of which can be found in Marde{\v{s}}ić's book, \cite{mardesic:strong:2000}. 
 \p
 The innovation in strong shape theory was, thus, to record the higher homotopies used in the approximation process. For instance, when using the \v{C}ech nerve of  open covers,  choices have to be made of the refinement maps between covers, and whilst these cannot be made to respect the iterated associativity required, `on the nose', so they do not give a commutative diagram indexed by the partially ordered set of open covers, they can be made \emph{coherently} with respect to further refinement, so recoding the higher homotopies in the approximation.  One replaces a homotopy approximation by a homotopy coherent one.
 \p
The link between that quite geometric topological approach and the notions emerging from Sugawara's work came from the work on homotopy everything algebraic structures by Boardman  and Vogt, \cite{boardmanvogt}, and then in Vogt, \cite{vogt:1973} which gave the link between a detailed `geometric' approach to homotopy coherent diagrams and a homotopical approach related to model category theory.  We will see links to Vogt's results later on.
  \end{rmk}
   
   

 The associativity being central to Sugawara's work on homotopy multiplicativity \cite{sugawara:hc} this study, it was not long before  maps of spaces with associative structures were  generalized to maps of topological categories and \emph{homotopy coherent functors}, which 
arise in other contexts involving topological  and simplicially enriched categories, \cite{cordier-porter:coherent}.

\begin{dfn} For topological  categories $\mathcal C$ and $\mathcal D$, a \emph{ functor up to strong 
homotopy} $F$, also known as a \emph{ homotopy coherent functor}, consists of maps $F_0: Ob\ \mathcal C\to Ob\ \mathcal D, F_1:Mor \to Mor $ and
maps $F_n:I^{n-1}\times {\mathcal C}_n \to {\mathcal D}_n$ such that 
  \begin{alignat}{2}
  F_1(&x\to y) : F_0(x)\to F_0(y) \notag \\
F_p(t_1, &\cdots ,t_{p-1}, c_1,\cdots ,c_n) = F_{p-1}(\cdots
,\hat{t}_i,\cdots
,c_ic_{i+1},\cdots ) &&\qquad\mbox{if } t_i=0  \notag\\
&= F_i(t_1, \cdots ,t_{i-1}, c_1,\cdots ,c_i)F_{p-i}(t_{i+1},\cdots
,c_{i+1},\cdots ,c_p) &&\qquad\mbox{if } t_i=1.  \notag
  \end{alignat}
  \end{dfn}

\p 
 Three particularly interesting examples of topological categories and their use are given by 
\begin{itemize}
\item the path space $B^I$, see section \ref{B^I};
\item the singular complex of $B$, for which see section \ref{Sing(M)}
\item[]\hspace{-1cm} and 
\item a `good'  open cover of $B$, which is in section \ref{good open cover}
\end{itemize}
\p
These also can, and will, be considered for smooth manifolds and smooth maps.

\section{Background/Motivation: ... and for dg-categories?}
Of course, there are complete analogs of these ideas for categories of (co)chain complexes and, more generally, dg categories.
We note that for  dg categories, that is categories enriched over some category of chain complexes,  a chain analog of cubes applies, so the $I^{n-1}$ in the above is replaced by a chain complex  representing it.

\p
For an associative algebra $A$, an $A$-module $M$ can be considered as a representation of $A$, a map 
$A\otimes M\to M$ or $A\to End(M).$  In a differential graded context, one again can consider representations up to homotopy of $A$ on $M$.
On the chain (dg) level, the corresponding notion is related to that of a \emph{twisting cochain}, the twisted tensor product differential as
 introduced in \cite{brown59}, for
modeling the chains on the total space of a principal fibration in terms of chains on the base and chains on the fiber.
\p
One can up the ante further by introducing $A_\infty$-spaces, $\Aoo$-algebras and modules and further $\Aoo$-categories,
but take care: do you want objects and morphisms to be $\Aoo$ in an appropriate sense or only the morphisms?
 It is worth pointing out that $A_\infty$-maps between strictly associative dg algebras were studied before \cite{jds:hahI,jds:hahII} by Sugawara as parameterized by cubes and called  \emph{strong homotopy muliplicative} maps.  $A_\infty$-maps of  
 $A_\infty$-spaces are parameterized by more complicated polyhedra.
\p
\dbend 
\p
 Warning! The strongly related notions of covariant derivative, connection, connection form and parallel transport have distinct definitions, 
but many authors use the names interchangeably in stating theorems while many proofs depend on a particular definition.

\section{Homotopy twisting cochains}
The analog of a fiber bundle, $F\to E\to B$, in terms of simplicial sets, $\mathcal F \to \mathcal E \to \mathcal B$, was formalized by
 Barratt, Gugenheim and Moore \cite{BGM}
using a \emph{twisting function} and called a \emph{twisted cartesian product}\footnote{They mention the concept had been around for some years. At the time they wrote the nomenclature was semi-simplicial complex or even complete semi-simplicial.}, then Ed Brown \cite{brown59} constructed an 
 algebraic chain model of a twisted cartesian product, 
  a \emph{twisted tensor product}  of chain complexes, $(C_*B\otimes C_*F,D)$,
where  the differential $D$ is $d_B \otimes 1 + 1 \otimes d_F$ twisted by adding a \emph{twisting term} $\tau: C_*B\to C_* Aut F$, that is
$$
D = d_B \otimes 1 + 1 \otimes d_F + \tau
$$
 where $\tau$ corresponds to  a representation of
  the cobar algebra $\Omega C_*B$ on $F.$ 
  Here $C_*(B)$ is regarded as a dg coalgebra and $C_* Aut F$ as a dg algebra.
 The defining relation that $\tau$ satisfies is
 $$d_F \tau + \tau d_B = \tau \cup \tau.$$
 In turn, Szczarba \cite{szczarba61} showed that for a twisted cartesian product base, $B$, with group $G$, there is 
there is a twisting cochain in $Hom(C_*(A), C_*(G)).$ Kadeishvili \cite{kad:twisting} studied twisting \emph{elements} in relation to $a\smile_1 a$ in a  homotopy Gerstenhaber algebra 
 \cite{GerVor95}. Quite recently for G-bundles, Franz \cite{franz:szcz} proved that the map $\Omega C_*B \to C_* G$ 
is a quasi-isomorphism of dg bi-algebras. That $\Omega C_*B$ is a dg bi-algebra follows from work 
 of Baues \cite{baues:doublebar} and Gerstenhaber-Voronov  \cite{GerVor95} using the homotopy Gerstenhaber structure of $C_*B$.
 \p
 For an algebra $A$, there is a \emph{universal} twisting  morphism $\tau_A:BA\to A$
such that for any twisting morphism $\tau:C \to A$,  there is a unique
morphism $f_\tau:C \to BA$ of coalgebras with $\tau=\tau_A\circ f_\tau$ \cite{HMS74}.
\p
The above can be generalized to  \shreps\footnote{If we write $D\tau$ for $d_F \tau + \tau d_B,$ we can see this as another manifestation of the Maurer-Cartan principle,
 which has subsumed the integrablity condition in deformation theory and others.}, leading to a notion of a \emph{homotopy twisting cochain} with values in an $A_\infty$-algebra.
 \p
 At the end of section \ref{RH}, twisting cochains will be revisited in relation to $\infty$-local systems. 
 \section{The Grothendieck construction and his `pursuit of stacks'}
 We mentioned the classical case of covering spaces and the equivalence between such and actions / representations of the fundamental group. This has a purely algebraic aspect as it mirrors the structures evident in Galois theory. 
 \p
 The construction of the covering space from the action is a simple example of a type of semi-direct product.  Both Ehresmann and Grothendieck looked at the relevance of this for fiber spaces and fibered categories as Grothendieck had introduced in his work.  The notion of a fibred category mimics both topological fibrations and the relationship between the category of all modules (over all rings) and the category of rings itself. 
 \p
 Looking at the analog of local sections in this categorical setting, one finds that a fibred category $p:\mathcal{E}\to \mathcal{B}$ corresponds to a pseudo-functor
 from the base category to the category of small categories, picking out the fibers and the actions linking them together. What is a pseudo-functor? It is a functor up to `natural equivalence of functors'  and thus up to the relevant notion of `homotopy' in this setting.  This, thus,  is another instance of a homotopy coherent representation, this time of the category $\mathcal{B}$ into the 2-category of small categories. Of special interest for us is the case of categories \emph{fibred in groupoids}. Here the pseudo-functor  will take values in the category of small groupoids\footnote{so of homotopy 1-types}.
 \p
 Grothendieck's construction started with a pseudo-functor and constructed a fibred category. In other words, it started with a homotopy coherent representation of the base category and built a fibred category from it. 
In his extensive discussion document \cite{groth:PS}, Grothendieck started exploring the higher homotopy analogues of fibred categories, thus analogues of simplicial fibrations.  The fibres were to be models of homotopy types, extending both  covering spaces,  categories fibred in groupoids, and more generally  homotopy $n$-types, for any value of $n$ including $\infty$. These were his $n$-stacks\footnote{There is an extensive literature on this, but here is not the place to explore this, although it does relate to homotopy coherent representations.}.
 
 
  Again homotopy coherence of the action, of the representation of the base homotopy type is key to building the fibred structure over that base. 
 
\section{The (topological) path space $B^I$}\label{B^I}
Returning to the topological case, the topological path space $B^I$ can be regarded as the space of morphisms of a topological category with $Ob = B$.
\p
A Hurewicz fiber space $F\to E\to B$ is a continuous mapping satisfying the \emph{homotopy lifting property} with respect to any space:
 Given a homotopy $\lambda:I\times X\to  B$ and an initial value $X\to F$, we have $\tau: I\times X \to E$ so that
\begin{displaymath}
\begin{CD}
   I\times   X @>\tau>>E\\
    @V{ }VV              @VV{p}V \\
I \times X  @> \lambda>>    B
\end{CD}
\end{displaymath}
\noindent is  commutative. 
\p
The action of the based loop space mentioned above 
extends (not uniquely) to an action  $Maps(I\to B) \times F_s\to F_t$ where for a path $\lambda: I\to B$ denote $F_s$ the fiber over $\lambda(0)$ and 
$F_t$ the fiber over $\lambda(1)$.
Non-tradtionally, this might be called \emph{parallel transport}.
\begin{dfn} \emph{Higher (parallel) transport} for a graded vector bundle $p:E\to M$ is a homotopy coherent functor $||:PM \to End(E)$.
\end{dfn}

 Arias Abad and  Schaetz \cite{as:parallel} use ``higher parallel transport'' to refer to  the \shrep denoting the category of representations up to homotopy of the simplicial set $Simp(M)$ of smooth
singular chains on $M$.
 
\section{The (smooth) path space $M^I$}
 
For vector bundles rather than covering spaces, according to n-lab:
\begin{quote}
Apparently one of the oldest occurrences of the idea that a principal bundle with connection  over a connected base space 
may be reconstructed from its holonomies around all smooth loops (for any fixed base point)
\end{quote}
 appeared in or was implied 
 by Kobayashi in 1954.
 
 
 \p
   The covering space  classification can be generalized  to give  an equivalence  between representations of  $\pi_1(M)$ 
   and  vector bundles $V\to E\to M$ with flat connection on $M$.
    A flat connection determines the lifting $\tau: I \to E$ uniquely, then
the \emph{holonomy}
with respect to a curve  is given by the evaluation of $\tau$ on  the path in $M$.
The \emph{holonomy group} is the image $\tau_\ast:\pi_1(M)\to GL(V)$   as a subgroup of the structure group of the bundle. 
The other direction, from  a representation  $\pi_1(M)\to GL(V)$ to a bundle with flat connection, 
is achieved by the associated $GL(V)$-bundle construction.
 \section{The  simplicial set/$\infty$-category $Simp(M)$ of $M$}\label{Sing(M)}
Just as classical holonomy is given by the evaluation of the lifting $\tau$, 
 consider  \emph{higher holonomy} by lifting  
$\sigma:\Delta^n\to M$ to $\tau^n:\Delta^n\times V\to E$ and related \emph{higher transport,}
  but now with  $E \to M$ a graded vector bundle with fiber $V$,  a  differential graded vector space.
  \p
  The idea of  higher holonomy was introduced by Chen as \emph{generalized holonomy\footnote{Only in \cite{Chen:iterated} does he use the word holonomy.}},
  elaborated   by  Igusa \cite{igusa:iterated} and  later  related to a notion of \emph{representation up to homotopy} \cite{as:parallel}.
 \p
  The technology of Chen's iterated integrals suggests a precise and rather natural way to handle \emph{ generalized holonomy.}  
 To capture such `higher structure', Chen used maps of  simplices $\sigma:\Delta^k\to M$.
 Denote the standard ordered $n$-simplex $\Delta^n$ as  $<0,1,\cdots,n>$ with vertices labelled $0,1,\cdots,n$.  Sub-simplices are denoted
 $<i_0,i_1,\cdots, i_j>$.  
 The face and degeneracy maps for a simplicial set  are:
\begin{align*}
 \partial_q <0,1,\cdots,n> &= <0,1,\cdots,q-1,q+1,\cdots, n>\\
 s_q<0,1,\cdots,n>& =<0,1,\cdots,q,q,\cdots,n>
\end{align*}
 \begin{dfn} $Simp(M)$ is  the simplicial set of (smooth) maps of  simplices $\sigma:\Delta^k\to M$ 
  For   $\sigma:\Delta^k\to M$, we denote by $V_i$ the fibre over the image of the vertex $i\in \Delta$.
  \end{dfn}
 \begin{dfn}\label{pathsinDelta} A \emph{representation up to homotopy}\footnote{There is a similar but earlier use of that  name \cite{jds:BAMS}.
Notice the `cubical' nature of the condition.} of $Simp(M)$ on a graded vector space $V$  is a collection $\theta$ of maps
$\{ \theta_{k}\}_{k\geq 0}$ which assign to any $k$-simplex $\sigma:\Delta^k\to M$ a map $ \theta_{k}(\sigma): I^{k-1}\times V \to V$ satisfying, for any $v\in V$, the relations:
\vskip1ex
\noindent
$\theta_0$ is the identity on  $V$ 
\vskip1ex
\noindent $\theta_k(\sigma)(t_1,\cdots,t_{k-1}, -): V \to V$ is an isomorphism for any $(t_1,\cdots,t_{k-1})\in I^{k-1}$.
\vskip1ex
\noindent For any $1\leq p\leq k-1$ and $v\in V$, 
\noindent
$$ \theta_{k}(\sigma)(\cdots ,  t_p=0, \cdots, v) =  \theta_{k-1}(\partial_p\sigma)(\cdots,\hat t_p,\cdots,, v)$$
$$ \theta_{k}(\sigma)(\cdots ,  t_p=1, \cdots, v) = $$ 
$$\theta_{p}(<0,\cdots,p>)\Big(t_1,\cdots, t_{p-1},  \theta_{q}(<p,\cdots,k>)(t_{p+1},\cdots,t_k,v)\Big).$$

\end{dfn}

This collection of maps is \emph{coherent} in the sense  that it respects the facial structure of the cubes and of the simplices.
\begin{theorem}For any flat graded vector bundle $p:E\to M$ with a flat  graded connection, there is a representation up to homotopy $\theta$ of $Simp(M)$ on $E$. 
\end{theorem}
 The desired maps $\theta_n$  are constructed realizing any simplex $\Delta^n$ as a family of paths with fixed endpoints $0$ and $n$: 
$\gamma_n:I^{n-1}\to P\Delta^n$ where
\noindent $\gamma_1(0)$ is the trivial path, constant at $0$
and $\gamma_2:I \to \Delta^1$ is the `identity'. 

\noindent For any $1\leq p\leq k-1$, 
$$\gamma_{k}(\cdots ,  t_p=0, \cdots) = \gamma_{k-1}(\cdots,\hat t_p,\cdots)$$
\noindent and
$$\gamma_{k}(\cdots ,  t_p=0, \cdots) = 
\gamma_p(t_1,\cdots, t_{p-1})\gamma_{q}(t_{p+1},\cdots,t_{k-1}).$$

Such maps were first produced by Adams \cite{adams:cobar} in the topological context by induction using the contractability of $\Delta^n$.  Later specific formulas were introduced by Chen \cite{Chen:iterated} and then equivalently but more  transparently, by Igusa \cite{igusa:twisting}.

\begin{rmk}

Given a category $\Cal$, Leitch \cite{leitch:homotopy-cube}  constructs a category  with the same objects, but  for each morphism $\phi$ of $\Cal$,  there is a space of morphisms, which he called the derived space of $\phi$. Considering  $\Delta^n$ as a category,  the derived space of $\Delta^n$ turns out to be $I^{n-1}$. This was one of the precursors of Cordier's homotopy coherent nerve construction \cite{cordier82}, for which see section \ref{homotopy coherent nerve}.
\end{rmk}
Just as representations of  $\pi_1(M)$ are related to
 vector bundles with flat connection on $M$, there is an alternative way to look at homotopy representations
 in terms of differential forms on $PM$.
   Consider a graded vector bundle $p:E= \prod E^k \to M$.
Let  $End^p(E)$ denote the degree $p$ part of the endomorphism bundle of $E$:
$$End^p(E)=\prod\Hom(E^k,E^{k+p}).$$
A $\Z$-graded connection  \emph{form} is just the analog of a classical connection form, 
but with careful attention to grading and signs. 
\begin{definition}
A $\Z$-\emph{connection form} $A$ on a graded vector bundle $p:E\to M$ is a form of total degree $1,$  i.e.  in $\oplus \Omega^p(M;End^{1-p}(E))$.
\end{definition}
As such, it corresponds to (a family of) differential forms with values in $End(E)$.
Let $\Omega^\bullet(M)$  be
the graded algebra of smooth differential forms on $M$ and let $\Omega^\bullet(M;E)$
be the graded $\Omega^\bullet(M)$-module of forms with values in $End\  E$.
As usual, it is useful to describe graded connections
locally, so the appropriate covariant derivative can be written:
$ d + A$, where $A = A_0+A_1+A_2+\cdots$ with $A_p \in  \Omega^{p}(M;E^{1-p})$

For every $t\in I$, consider the evaluation map
\[
	ev_t:PM \to M
\]
sending $\g$ to $\g(t)$.
Let $W_t$ be the pull back of $E$ to $PM $ along $ev_t$. That is,
\[
	(W_t)_\g=V_{\g(t)}.
\]
For $0\le s\le t\le 1$, let
$\Hom^q(W_s,W_t)$
be 
the space of degree $q$ graded homomorphisms from $W_s$ to $W_t$. Define a smooth bundle
\[
	\Omega^p(PM ,\Hom^q(W_s,W_t))=\Omega^p(PM)\otimes\Hom^q(W_s,W_t)
\]
whose fiber over $\g$ is the vector space of smooth $p$-forms 
with coefficients in 
\[
	\Hom^q(W_s,W_t)_\g=\Hom^q(V_{\g(s)},V_{\g(t)}).
\]
 \begin{dfn}A family of \emph{homotopy coherent forms} on $PM$ means a family of forms $\Psi_p(s,t)$ for all $0\le s\le t\le 1$
  for  the smooth singular simplicial set $Simp(M)$ of $M$ such that 
\[
	\Psi_p(s,t)\in \Omega^p(PM ,\Hom^{-p}(W_s,W_t))
\]
satisfying the following at each $\g\in PM $:
\begin{enumerate}
\item $\Psi_0(s,s)_\g$ is the identity map in
\[
\Hom^0(V_{\g(s)},V_{\g(s)})
\]
\item $\Phi(s,t)_\g$ satisfies a certain  first order linear differential equation \cite{igusa:twisting}:
\[
	\Psi_p(\g,t,s)\in\Omega^p(PM,\Hom^{-p}(W_s,W_t)).
\]
\end{enumerate}
\end{dfn}
  \begin{thm}  \cite{igusa:twisting}:
Given a graded vector bundle $p:E\to M$  with a \emph{flat} $\Z$-\emph{connection form} $A$, there is a family of 
 \emph{homotopy coherent forms} on $V$  of the smooth singular simplicial set $Simp(M)$ of $M$
 \end{thm}
 The differential equation determines $\Psi_p(s,t)$ uniquely as a $p$-form on $PM$. 
Starting with $\Psi_0(s,t),$
 one can find $\Psi_p(s,t)$ by induction on $p$ using a version of Chen's iterated integrals.

\section{A `good' open cover $\mathcal U$}\label{good open cover}
An open cover  $\mathcal U =\{U_\alpha\}$ is `good' if all the $U_\alpha$ as well as all their intersections are contractible.
 If $\{ U\da\}$ is the open covering, the disjoint union
$\coprod U\da$ can be given a rather innocuous structure of a topological category $U$,
i.e., $Ob U = \coprod U\da$ and $Mor U=\coprod U\da\cap U\db$
that is $x\circ y=x=y$ is defined iff $x\in U\da$, $y\in U\db$ and $x=y$.  For an ordinary vector bundle
with structure group $G\subset GL(V)$, we have
 {\em transition functions} $g\dab : U\da\cap U\db \to G$ which satisfy the cocycle condition:
  $g\dab  g\dbgam =g_{\alpha\gamma}$,
 but for fiber spaces life is not so straightforward. The transition functions are defined in terms of 
 trivializations.
 \p 
 A differential graded vector bundle  $p: E\to M$ with fiber $V$, a differential graded vector space, is locally  trivial over some
open covering $\{ U\da\}$  with trivializations: dg maps $h\da$ and $k\da$ such that 
 $$h_\alpha:p^{-1}(U\da ) \rightleftarrows U\da\times V  : k\da$$
 where $h\da$ and $k\da$ are inverse fibrewise equivalences.
\p   The transition functions are again defined by the equation
  \[
h\da k\db (x,v) = \bigl( x, g\dab (x)(v)\bigr) , \qquad x\in U\da\cap U\db ,\; v\in V,  \]
but now it is not necessarily true that $g\dab  g\dbgam =g_{\alpha\gamma}$; rather
instead of the cocycle
condition, one obtains only that $g\dab g\dbgam$ is homotopic to $g\dagam$ as a map on
$U\da\cap
U\db\cap U\dgam$.  Moreover, on multiple intersections, higher homotopies arise and constitute  a {\em homotopy coherent functor},
relevant to classifying the fibration.
\p
The collection of higher homotopies was called a \emph{homotopy transition cocycle} \cite{wirth:diss,jw-jds}.
They determine a fibration (up to appropriate equivalence) by assembling pieces $U\times V$, not by glueing them directly to each other, but only after inserting `connective' tissue between those pieces , e.g. consider ,
$$U\da \times V \cup (I\times U\da \cap U\db) \cup U\db\times V.$$
\section{Variations on the theme of a homotopy coherent nerve}\label{homotopy coherent nerve}
A useful tool in handling the properties of homotopy coherent functors is to view the various forms of the homotopy coherent nerve construction.
\p
The definition in probably the most standard form is for a simplicial enriched category, $\mathcal{C}$, in which the hom-sets between the objects are Kan complexes. This fits the cases of $\mathcal{C}=Top$ or $Kan$ and can also  be adapted to handle the category of chain complexes with a bit of extra work.  The basic construction is due to Cordier \cite{cordier82}, using ideas developed by Boardman and Vogt, \cite{boardmanvogt}.
\p
The idea is that for each ordinal $[n]$, one builds a simplicial enriched category, $S[n]$, that in a certain sense resolves the category $[n]$, and such that $\mathcal{S}-Cat(S[n], \mathcal{C})$ gives the collection of all homotopy coherent diagrams within $\mathcal{C}$ having the form of an $n$-simplex. These collections form a simplicial class, (they will often be too large to be a simplicial set), called the homotopy coherent nerve of the $\mathcal{S}$-category $\mathcal{C}$ and denoted $Ner_{h.c.}(\mathcal{C})$. It then is clear that a homotopy coherent diagram having the  form of a small category $\mathsf{A}$ is given by a simplicial map from $Ner(\mathsf{A}$ to $Ner_{h.c.}(\mathcal{C})$.
\p
Vogt proved, \cite{vogt:1973}, for $\mathcal{C}=Top$, that there was a category of homotopy coherent diagrams of type $\mathsf{A}$ in $Top$\footnote{The definition and composition of homotopy coherent morphisms is difficult to describe briefly.} and that this was equivalent to the category obtained from $\mathcal{C}^\mathsf{A}$ by inverting the `pointwise' homotopy equivalences. This gives a number of useful links between the explicit definition of homotopy coherent maps and their interpretation in more homotopical terms.
\p
This homotopy coherent nerve is a quasi-category, so is a model for an $\infty$-category.  Given any quasi-category $K$, and thinking of this object as an $\infty$-category, we define a homotopy coherent functor from $K$ to $\mathcal{C}$ to be simply a simplicial map from $K$ to $Ner_{h.c.}(\mathcal{C})$. In particular for $K= Sing(X)$, which, as it is a Kan complex is a quasi-category, and is known as the $\infty$-groupoid of $X$, a $\infty$-local system in $\mathcal{C}$ is simply a mapping from $Sing(X)$ to $Ner_{h.c.}(X)$. 
These constructions can be adapted for dg-categories, as was done by Lurie, \cite{lurie:higher:2011}. It can further be adapted to give a homotopy coherent nerve of an $A_\infty$-category as in Faonte \cite{faonte:simplicial-nerve:2013}, and we meet a variant of this in the next section.

\section{$\mathsf{Loc}^{\mathcal C}(\pi_{\infty}M)$}\label{RH}
In a  variant of the classical Riemann--Hilbert equivalence, a map
$$ Flat(M) \rightarrow Reps(\pi_1(M))$$
is developed by calculating the holonomy with respect to a flat connection.  The holonomy descends to a representation 
of $\pi_1(M)$ as a result of the flatness. 
  In their generalization of the  Riemann--Hilbert correspondence, \BS  write:
\begin{quote}
 From the perspective of (smooth) homotopy theory,  the manifold $M$ can be replaced by its infinity-groupoid $\pi_{\infty}M := Sing^{\infty}M $ of smooth simplices.  Considering the correct notion of a representation of $\pi_{\infty}M$ will allow us to produce an untruncated Riemann-Hilbert theory.  
 \end{quote}
The notion they define is  that of an
\emph{infinity-local system} on $M$, denoted $\mathsf{Loc}^{\mathcal C}(\pi_{\infty}M),$
which is essentially a homotopy coherent functor from the $\infty$-category $Simp(M)$ to the category of 
${\mathcal Ch}$ of chain complexes.
\begin{rmk}
Rivera and Zeinalian \cite{RZ:ttp} use both $\infty$-local systems and the realization of a simplex $\Delta^n$ as a family of paths 
parameterized by $I^{n-1}$ (see \ref{pathsinDelta}) to shed new light on Brown's classical result \cite{brown59}.
\end{rmk}
To generalize the classical  \RH equivalence, \BS  invoke   the Serre-Swan theorem to change the category of flat graded vector bundles to that of ``perfect modules with flat $\mathbb Z$-connection. They then write:
\begin{quote}
It would be an interesting problem in its own right to define an inverse functor which makes
use of a kind of associated bundle construction. 
\end{quote}
 That is work for the future.

\bibliographystyle{amsplain}

\bibliography{Loriano-hcf4}
\ddo
 \p
 
The classical  \RH equivalence 
involves the functor
\begin{equation*}
 Reps(\pi_1(M)) \rightarrow Flat(M).
\end{equation*}
which is achieved by the associated bundle construction for the structure group $GL(V)$.

\p
Block and Smith work in the category of  

but remark:
\begin{quote}
It would be an interesting problem in its own right to define an inverse functor which makes use of a kind of associated bundle construction. 
\end{quote} This is exactly what we do next, following the unpublished thesis of Surapaneni \cite{Surapanenini}.

\section{Constructing a graded bundle with a flat connection}

\trd{Above we have tried to act as translator of terminologies that in fact refer tgo equivalent concpts}


\begin{thm}
Given a representation up to homotopy of $Simp(M)$ on ???? $V\times M$

 there is a graded  bundle $E$ with a flat connection 
for which the representation up to homotopy constructed as in Theorem
is equivalent to the given one.\end{thm}

\trd{OR the family form version?? OR the || version?}

  Given a dg vector space $(V,d)$ and $\sigma:<1,\cdots,n> \to M$, let $E_\sigma$ be the pullback of the trivial bundle $V\times M$.
 \p
 Alternative: let $E_\sigma$ be the trivial bundle $V\times \sigma$.
\p
Alternative 2: Surapenini version
\p
\ddo
These pieces don't fit together well, so we use the data from $\theta$  to attach additional connective `tissue' to  the disjoint union $\coprod E_\sigma$
 to form a bundle $E_{Simp}$ 
 over the geometric realization of $Simp(M)$.
 \p
   Recall the  mapping cylinder construction and apply it  sequentially as follows:
  For $n=0$, let
 $E_0 = \coprod E_i =\coprod V_i$:
   For each 1-simplex $\sigma = <i i+1>$, \trd{0 1 sufficient ?} attach  
   the mapping cylinder of $\theta_1(\sigma)$ from $V_i$ to $V_{i+1}$.
   
   \trd{sketch}
 \p
 Then,   for each 2-simplex $\sigma=<012>$, attach  $E_\sigma\times I^2$ as the union of the mapping cylinders of $\theta_2$ for each face of $I^2$.
 Of course, we must check this is well defined at the corners of $I^2$.

OR $E_0 = V_i$ and $E_1 = I\times V/~$


 OR
 the mapping cylinder of Igusa's / to  
 Boundary of $I^2$ in terms of $\theta$'s maps through $\Delta^2$?
 
 02 ----- 12 (01) given by 012 of Igusa
 
 V---------01
 
 See 10 AM 12/24 and circ(2) for generic doub;le mapping
 
 Let Hom(E,E) be the vector bundle over $M\times M$ such that the fiber over $(x,y)$ is $Hom (E_x,E_y)$
 
 Denote by $P^2E \to PM$ the subbundle of $PE$ which is the pullback of $PM\to M\times M$ as for a category.

associated transport $\tau:PM\to P^2M$

\ddo